\documentclass[10pt]{article}
\usepackage{amsthm,amsmath,amssymb,amscd,verbatim,epsfig}
\newtheorem{thm}{Theorem}
\newtheorem{lem}{Lemma}

\newtheorem{defn}{Definition}

\newtheorem{rem}{Remark}

\begin{document}
\title{A refinement of 
Sharkovskii's theorem on orbit types characterized by two parameters}

\author{Bau-Sen Du\\
Institute of Mathematics \\
    Academia  Sinica \\
    Taipei 115, TAIWAN \\ 
    dubs@math.sinica.edu.tw \\
	and \\
	Ming-Chia Li\\
	Department of Mathematics\\ National Changhua University of Education \\
	Changhua 500, TAIWAN \\
	mcli@math.ncue.edu.tw\footnote{current address: Department of Applied Mathematics, National Chiao Tung 
University, Hsinchu 300, TAIWAN, mcli@math.nctu.edu.tw}
}

\date{}

\maketitle

\begin{abstract}
The so-called type problem or forcing problem is considered as a way to generalize Sharkovskii's theorem.  In this 
paper, by focusing on certain types of orbits, we obtain a solution of the type problem, which gives a refinement of 
Sharkovskii's theorem on orbit types characterized by two parameters. 
\end{abstract}


\section{Introduction}
For continuous interval maps, in 1964, the remarkable theorem of Sharkovskii \cite{Sharkovskii1964} gave the 
complete answer of the following question: {\em Given a periodic orbit of a specified period, find the other periods 
of periodic orbits that must exist} (Theorem~\ref{thm1}).    
One can classify types of orbits, up to symmetry, depending on the arithmetic ordering of the points in the real 
line; refer to \cite{Baldwin1987}.  For example, two  period-$4$ orbits $\{x_0\to x_1\to x_2\to x_3\to x_0\}$, one 
with $x_0<x_1<x_2<x_3$ and the other one with $x_0<x_2<x_1<x_3$, are considered to have different types.  
The so-called {\em type problem} or {\em forcing problem} is the following: 
\begin{enumerate}
\item[]{\em Given a period-$n$ orbit of a specified type, find, for any positive integer $m$, the types of 
period-$m$ orbits that must exist.}
\end{enumerate}
There are at least ${(n-1)!}/2$ different types of period-$n$ orbits.  Due to the factorial  growth of the number of 
different types, the complete solution of the type problem  
is not an easy task.  The type problem is closely related to bifurcations of any one-parameter family of continuous 
interval maps in the sense that its solutions can show all possible routes of bifurcations of various types of 
periodic orbits.  For more discussions on the type problem, refer to \cite{AlsedaLlibreMisiurewicz2000} and 
\cite{BlockCoppel1992}.

By considering the implications between existence of different orbit types characterized by one parameter, a well 
known solution of the type problem gives a refinement of Sharkovskii's theorem (Theorem~\ref{thm2}).  In this paper, 
we extend the result by considering different orbit types characterized by two parameters (Theorem~\ref{thm3}).

\section{Definitions and Statements of Theorems}
Let $I$ be a nontrivial compact interval and let $f$ be a continuous map from $I$ into itself. 

First of all, we give some basic definitions. The {\it forward orbit} of $x_0$ for $f$ is defined to be the sequence 
$\{x_i: i\geq 0\}$, where $x_i=f^i(x_0)$, i.e., the $i$-th iterate of $f$ at $x_0$.  A sequence $\{x_{-i}: i\geq 
0\}$ is called a {\em backward orbit} of $x_0$ for $f$ if $f(x_{-i})=x_{-i+1}$ for all $i\geq 1$.  It is possible 
that a point has several backward orbits if  $f$ is not one-to-one.  
A point $x_0$ is called a {\em period-$n$ point} of $f$ if its forward orbit $\{x_i\}$ 
satisfies $x_n=x_0$ and $x_i\neq x_0$ for all $0<i<n$.  A period-$1$ point is also called a {\em fixed point}.  

We say that the {\em property $P(n)$ holds} if $f$ has a period-$n$ point, and denote that {\em $P(n)\to P(m)$} if 
the property $P(n)$ implies the property $P(m)$. {\em Sharkovskii's theorem} reads as follows. 

\begin{thm}[\cite{Sharkovskii1964}]\label{thm1} Let $f$ be a continuous map from $I$ into itself.  Then the 
following diagram holds: 
\begin{align*}
&P(3)\to P(5)\to P(7)\to \cdots \\
\to &P(2\cdot 3)\to P(2\cdot 5)\to P(2\cdot 7)\to \cdots \\
\to &P(2^2\cdot 3)\to P(2^2\cdot 5)\to P(2^2\cdot 7)\to \cdots \\
\to &\cdots\\
\to &\cdots\to P(2^2)\to P(2)\to P(1).
\end{align*}
\end{thm}

In the following, some properties for $f$ with specific types of orbits are defined. 
\begin{defn}
Let $k$ and $n$ be positive integers.  We say that, for $f$, 
\begin{enumerate}
\item the property $L^k(n)$ holds if $f^k$ has a period-$n$ point $x_0$ with the forward orbit $\{x_i\}$ of the type 
either 
\begin{align*}
\text{$x_0<x_1<x_2<\cdots<x_{n-1}$}
\end{align*}
 or all inequalities reversed;
\item the property $S^k(2n+1)$ holds if $f^k$ has a period-$(2n+1)$ point $x_0$ with the forward orbit $\{x_i\}$ of 
the type either 
\begin{align*}
\text{$x_{2n}<x_{2n-2}<\cdots<x_2<x_0<x_1<x_3<\cdots<x_{2n-3}<x_{2n-1}$}
\end{align*}
 or all inequalities reversed; and 
\item the property $L^k(\infty)$ holds if $f^k$ has a fixed point $x_0$ with a backward orbit $\{x_{-i}\}$ of the 
type either 
\begin{align*}
\text{$x_0<\cdots<x_{-i}<x_{-i+1}<\cdots<x_{-2}<x_{-1}$}
\end{align*} 
 or all inequalities reversed.
\end{enumerate}
\end{defn}

In \cite{BhatiaEgerland1988}, it is proved that $L^k(\infty)\to L^k(n)$ and $L^k(n+1)\to L^k(n)$; see also 
\cite{Carvalho1989} for independent work.  In \cite{BlockCoppel1992}, it is shown that $S^k(2n+1)\to S^k(2n+3)$.  
Combining these results, together with Sharkovskii's theorem \cite{Sharkovskii1964} (see also \cite{Stefan1977}), 
one easily has the following refinement of Sharkovskii's theorem. 

\begin{thm}
[\cite{BhatiaEgerland1988},\cite{BlockCoppel1992},\cite{Sharkovskii1964},\cite{Stefan1
977}]\label{thm2} Let $f$ be a continuous map from $I$ into itself.  Then the following diagram holds: 
\begin{align*}
&L^1(\infty)\to \cdots\to L^1(5)\to L^1(4)\to L^1(3) \\
\leftrightarrow &P(3)\leftrightarrow S(3)\to P(5) \to S(5)\to P(7)\to S(7)\to \cdots \\
\to &L^2(\infty)\to \cdots\to L^2(5)\to L^2(4)\to L^2(3)\\
\leftrightarrow &P(2\cdot 3)\leftrightarrow S^2(3)\to P(2\cdot 5) \to S^2(5)\to P(2\cdot 7)\to S^2(7)\to \cdots  \\
\to &L^{2^2}(\infty)\to \cdots\to L^{2^2}(5)\to L^{2^2}(4)\to L^{2^2}(3) \\
\leftrightarrow &P(2^2\cdot 3)\leftrightarrow S^{2^2}(3)\to P(2^2\cdot 5) \to S^{2^2}(5)\to P(2^2\cdot 7)\to 
S^{2^2}(7)\to \cdots \\
\to &\cdots \\
\to &\cdots\to P(2^2)\to P(2)\to P(1).
\end{align*}
\end{thm}

We consider more properties for $f$ with certain types of orbits.  

\begin{defn}
Let $k$, $m$ and $n$ be positive integers.  We say that, for $f$, 
\begin{enumerate}
\item the property $L^k(m, n)$ holds if $f^k$ has a period-$(m+n)$ point $x_0$ with the forward orbit $\{x_i\}$ of 
the type either 
\begin{align*}
\text{$x_{m+n-1}<\cdots<x_{n+1}<x_n<x_0<x_1<x_2<\cdots<x_{n-1}$}
\end{align*}
 or all inequalities reversed;
\item the property $L^k(m,\infty)$ holds if $f^k$ has a fixed point $x_0$ with a backward orbit $\{x_{-i}\}$ of the 
type either 
\begin{align*}
\text{$x_{-1}<x_{-2}<\cdots<x_{-m}<x_0<\cdots<x_{-i}<x_{-i+1}<\cdots<x_{-m-1}$}
\end{align*}
 or all inequalities reversed; and 
\item the property $L^k(\infty,\infty)$ holds if $f^k$ has a fixed point $x_0=y_0$ with two backward orbits 
$\{x_{-i}\}$ and $\{y_{-i}\}$ of the type either 
\begin{align*}
\text{$x_{-1}<\cdots<x_{-i+1}<x_{-i}<\cdots<x_0=y_0<\cdots<y_{-i}<y_{-i+1}<\cdots<y_{-1}$}
\end{align*} 
 or all inequalities reversed.
\end{enumerate}
\end{defn}

Note that the properties $L^k(m, n)$ and $L^k(n, m)$ are equivalent. 

Now we state the main result of this paper. 
\begin{thm}\label{thm3}Let $f$ be a continuous map from $I$ into itself.  Then for all integers $k\geq 1$ and $i\geq 
0$, the following two diagrams hold: \\

\centerline{\epsfig{file=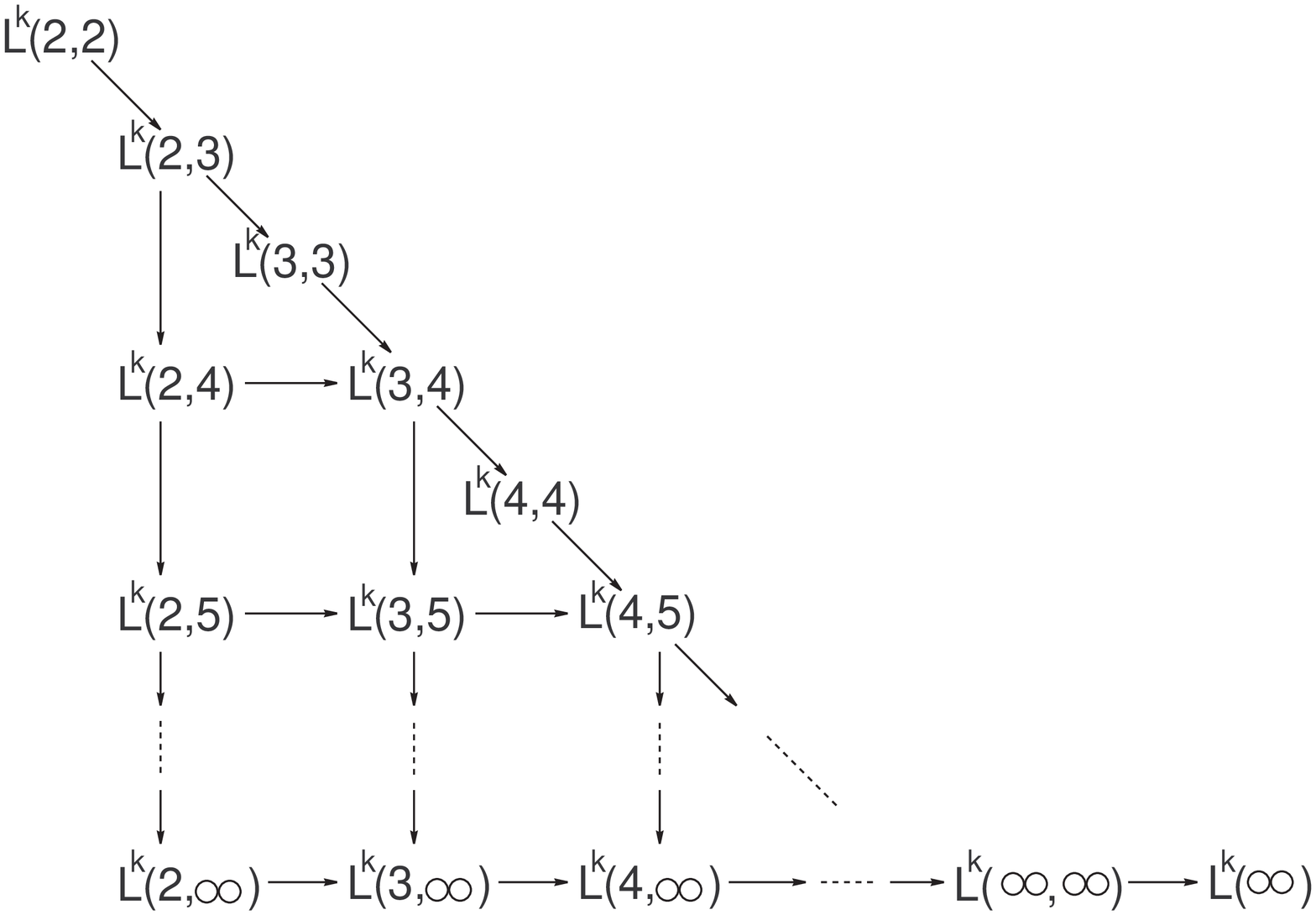,width=9.7cm,height=5.7cm}}

and 

\hspace{1cm}\centerline{\epsfig{file=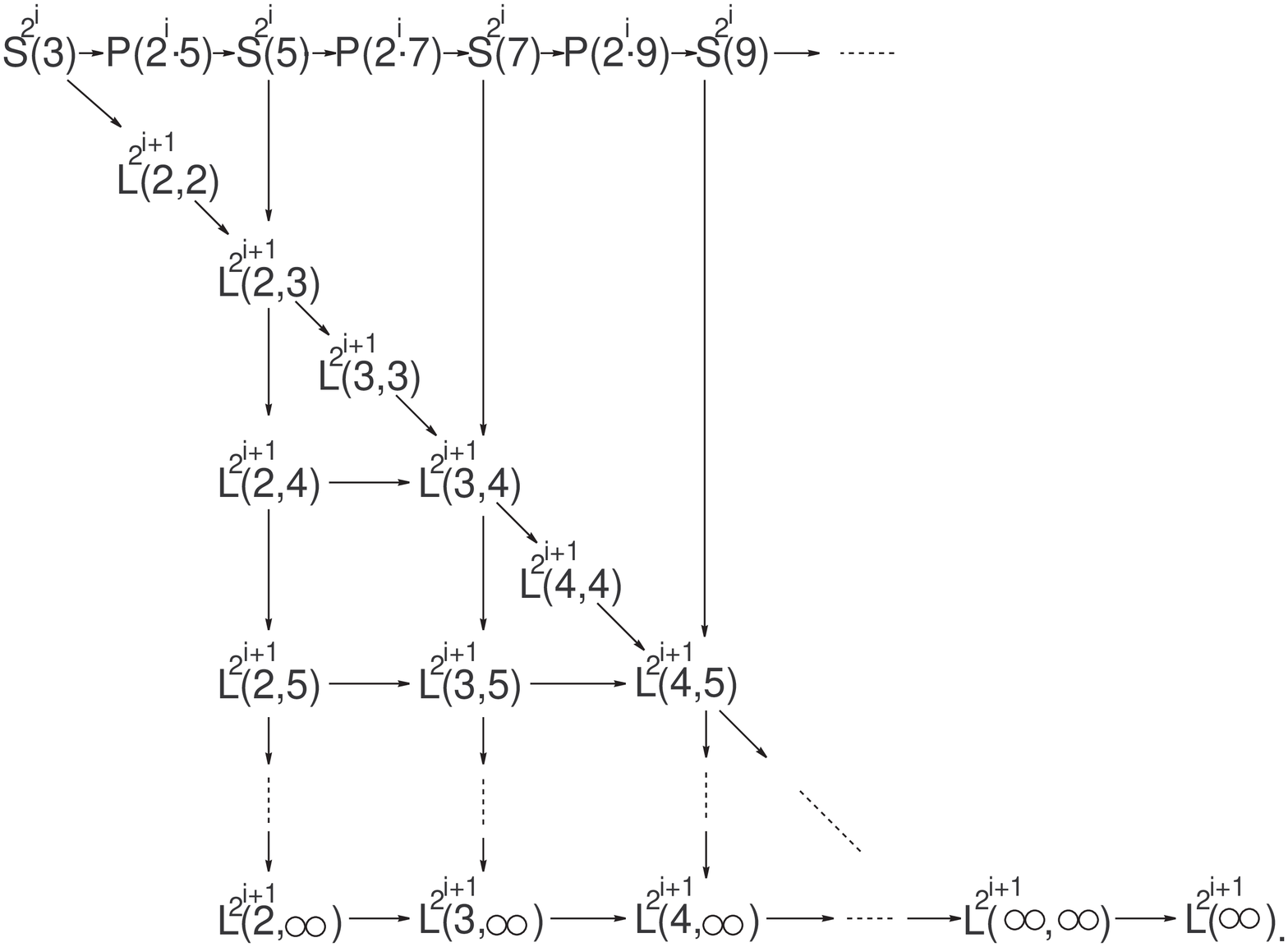,width=12.7cm,height=5.7cm}}
\end{thm}

\begin{rem}The above theorem, together with Theorem~\ref{thm2}, gives a refinement of Sharkovskii's theorem on orbit 
types characterized by two parameters.
\end{rem}

\section{Proof of Theorem~\ref{thm3}}
First, we recall a basic lemma, the proof of which can be found in \cite{Robinson1999}. 
\begin{lem}\label{lemma}Let $f$ be a continuous map from $I$ into itself.  Then the following statements are true. 
\begin{enumerate}
\item[$(a)$] If $J$ is a 
closed subinterval of $I$ with $f(J) \supset J$, then there exists a 
fixed point of $f$ in $J$.
\item[$(b)$] Let $J_i$, $0\leq i\leq n-1$, be closed subintervals of $I$.  If $f(J_{i})\supset J_{i+1}$ for all 
$0\leq i\leq n-2$  and $f(J_{n-1})\supset J_0$, then there exists a periodic point 
$y$ of $f$ in $J_0$ such that $f^i(y) \in J_i$ for all $1 \leq i \leq
n-1$ and $f^n(y) = y$.  
\end{enumerate}
\end{lem}

\medskip

The statements $S^{2^i}(2m+1)\to P(2^i\cdot (2m+3))\to S^{2^i}(2m+3)$ for all integers $i\geq 0$ and $m\geq 1$ come 
trivially from Theorem~\ref{thm2}. 

To prove the theorem, we only need to show that, for all integers  $m\geq 2$, $n\geq 2$, $k\geq 1$ and $i\geq 0$, 
\begin{enumerate}
\item[(i)] $L^k(m,n) \to L^k(m,n+1)$ and $L^k(m,n) \to L^k(m+1,n)$; 
\item[(ii)]$L^k(m,n) \to L^k(m,\infty) \to L^k(m+1,\infty) \to L^k(\infty,\infty) \to 
L^k(\infty)$; 
\item[(iii)]$S^{2^i}(3) \to L^{2^{i+1}}(2,2)$; and 
\item[(iv)]$S^{2^i}(2n+1) \to L^{2^{i+1}}(n,n+1)$.
\end{enumerate}

It suffices to prove items (i) and (ii) with $k=1$ and items (iii) and (iv) with $i=0$ because the other cases with 
$k>1$ and $i>0$ follow immediately by considering $f^k$ and $f^{2^i}$, respectively.

For item (i) with $k=1$, we may assume that $f$ has a period-$(m+n)$ point $x_0$ with the forward orbit $\{x_i\}$ of 
the type $x_{m+n-1}<\cdots<x_{n+1}<x_n<x_0<x_1<x_2<\cdots<x_{n-1}$.  Since $f([x_n, x_0]) \supset [x_n, x_0]$, 
Lemma~\ref{lemma}(a) implies that $f$ has a fixed point $z$ in $[x_n, x_0]$.  In fact, $z$ lies in $(x_n,x_0)$ since 
$x_0$ and $z$ have different periods.  Since $f([z,x_0])\ni x_0$, there is a point $w \in (z,x_0)$ such
that $f(w) = x_0$.  Now let $J_0 = [z,w]$, $J_1 = [w,x_0]$, 
$J_i = [x_{i-2}, x_{i-1}]$ for $2\leq i \leq n$, $J_{n+1} = [x_n, z]$, and $J_i = [x_{i-1}, x_{i-2}]$ for $n+2 \leq 
i \leq m+n$.  Then $f(J_i)\supset J_{i+1}$ for $0\leq i\leq m+n-1$ and $f(J_{m+n})\supset J_0$. By 
Lemma~\ref{lemma}(b), $f$ has a periodic point $y$ in $J_0$ such that $f^i(y)\in J_i$ for all $1\leq i\leq m+n$ and 
$f^{m+n+1}(y)=y$. It is clear that $y$ is neither $z$ nor $w$.  Since $f^i(y)\in J_i$ for all $0\leq i\leq m+n$, $y$ 
is a period-$(m+n+1)$ point and its forward orbit has the type so that the property $L^1(m,n+1)$ holds. This 
completes the proof of the first statement of item (i).  By the definition, the properties $L^1(m, n)$ and $L^1(n, 
m)$ are equivalent.  Thus the second statement follows from the first one.

For item (ii) with $k=1$, under the same assumption as above, we have shown the existence of a fixed point $z$ in 
$(x_n, x_0)$.  First, we find a backward orbit $\{z_{-i}\}$ of $z_0=z$ with the type so that the property $L^1(m, 
\infty)$ holds. 
Since $f^i((x_{m+n-i},x_{m+n-i-1}))\ni z_0$ for all $1\leq i\leq m-1$, we can, by induction on $i$, find  $z_{-i}$ 
with $f^i(z_{-i})=z_0$ in $ (x_{m+n-i}, x_{m+n-i-1})$ for all $1\leq i\leq m-1$.  Since $z_{-m+1}\in (x_{n+1}, 
x_n)\subset (x_{n+1}, z_0)\subset f((x_n, z_0))$, we can find $z_{-m}$ in $(x_n, z_0)$.  Again, by induction on $i$, 
we can find $z_{-i}$ in $(x_{m+n-i-1}, x_{m+n-i})$ for $m+1\leq i\leq m+n-1$, since $f^{i-m}(x_{m+n-i-1}, 
x_{m+n-i})\ni z_{-m}$.  Since $f((z_0,z_{-m-n+1}))\ni z_{-m-n+1}$, we can find $z_{-m-n}$ in $(z_0,z_{-m-n+1})$; 
continuing this process, by induction on $i$, we can find $z_{-i}$ in $(z_0, z_{-i+1})$ for all $i\geq m+n+1$.  We 
have proved that the property $L^1(m,\infty)$ holds.  Next, we find another backward orbit $\{z^\prime_{-i}\}$ of 
$z^\prime_0=z$ with type so that the property $L^1(m+1, \infty)$ holds.  Let $z^\prime_{-i}=z_{-i}$ for $1\leq i\leq 
m$. Since $f((z^\prime_{-m}, z^\prime_0))\ni z^\prime_{-m}$, we can find $z^\prime_{-m-1}$ in $(z^\prime_{-m}, 
z^\prime_0)$.  By induction on $i$, we can find $z^\prime_{-i}$ in $(z_{-i-1}, z_{-i})$ for all $i\geq m+1$, since 
$f((z_{-i-1}, z_{-i}))\ni z^\prime_{-i}$. This shows the truth of the property $L^1(m+1, \infty)$.  
By using the same argument, one can show that the property $L^1(\infty,\infty)$ holds.  By the definition, we have 
$L^1(\infty, \infty)\to L^1(\infty)$.  The proof of item (ii) is complete.    

For item (iii) with $i=0$, we may assume that $f$ has a period-$3$ point $x_0$ with the forward orbit $\{x_i\}$ of 
the type $x_2<x_0<x_1$.  Since $f((x_0, x_1))\ni x_0$, there is a point $w\in (x_0, x_1)$ such that $f(w)=x_0$.  Let 
$g=f^2$.  Then $g(x_2)=x_1$, $g(x_0)=x_2$, $g(w)=f(x_0)=x_1$, and $g(x_1)=x_0$.  By using Lemma~\ref{lemma}(a), 
there are fixed points $a, b$ and $c$ for $g$ such that $x_2<a<x_0<b<w<c<x_1$.  Let $J_0=[x_0, b]$, $J_1=[x_2, a]$, 
$J_2=[w,c]$, and $J_3=[c,x_1]$.  Then $g(J_i)\supset J_{i+1}$ for $0\leq i\leq 2$ and $g(J_3)\supset J_0$.  
Lemma~\ref{lemma}(b) implies that $g$ has a period-$4$ point in $J_0$ so that the property $L^2(2,2)$ holds.  The 
proof of item (iii) is complete. 

For item (iv) with $i=0$, we may assume that $f$ has a period-$(2n+1)$ point $x_0$ with the forward orbit $\{x_i\}$ 
of the type $x_{2n}<x_{2n-2}<\cdots<x_2<x_0<x_1<x_3<\cdots<x_{2n-3}<x_{2n-1}$.  It is clear that the property 
$L^2(n+1,n)$ holds, which is equivalent to the truth of the property $L^2(n, n+1)$, and so the proof of item (iv) is 
complete. 

We have finished the proof of Theorem~\ref{thm3}.

\noindent
{\bf Acknowledgment. } The authors would like to thank the referee for valuable suggestions which led to an 
improvement of this paper.

\bibliographystyle{amsplain}
\providecommand{\bysame}{\leavevmode\hbox to3em{\hrulefill}\thinspace}

\end{document}